\documentclass[12pt,reqno]{amsart}
\usepackage{amssymb}
\newtheorem{theorem}{Theorem}[section]
\newtheorem{proposition}{Proposition}[section]

\newtheorem{lemma}{Lemma}[section]
\newtheorem{definition}{Definition}[section]

\newfont{\bb}{msbm10 at 12pt}

\newcommand{\mysection}[1]{\section{#1}\setcounter{equation}{0}}

\def\pf{{\textit {Proof :} }}

\def\R{\hbox{\bb R}}

\def\P{\hbox{\bb P}}
\def\L{\hbox{\bb L}}

\def\N{\hbox{\bb N}}
\def\X{\hbox{\bb X}}

\def\i{\it i}

\def\E{\mathcal E}

\def\P{\mathcal P}
\newcommand{\bal}{\begin{align}}     \newcommand{\eal}{\end{align}}
\newcommand{\ba}{\begin{array}}      \newcommand{\ea}{\end{array}}
\newcommand{\bc}{\begin{center}}     \newcommand{\ec}{\end{center}}
\newcommand{\be}{\begin{enumerate}}  \newcommand{\ee}{\end{enumerate}}
\newcommand{\beq}{\begin{eqnarray}}  \newcommand{\eeq}{\end{eqnarray}}
\newcommand{\beQ}{\begin{eqnarray*}} \newcommand{\eeQ}{\end{eqnarray*}}
\newcommand{\bi}{\begin{itemize}}    \newcommand{\ei}{\end{itemize}}
\newcommand{\bt}{\begin{tabular}}    \newcommand{\et}{\end{tabular}}
\newcommand{\bdm}{\begin{displaymath}} \newcommand{\edm}{\end{displaymath}}

\newcommand{\rw}{\rightarrow}

\let\pa=\partial

\def\qed{\hfill{Q.E.D.}\smallskip}

\begin{document}

\title[Positivity of the Bondi mass]
{Positivity of the Bondi mass in Bondi's radiating spacetimes}

\author{Wen-ling Huang}
\address[Huang]{Department Mathematik, Schwerpunkt GD, Universit\"at Hamburg,
Bundesstr. 55, D-20146 Hamburg, Germany}
\email{huang@math.uni-hamburg.de}
\author{Shing Tung Yau}
\address[Yau]{Department of Mathematics, Harvard University,
Cambridge, MA 02138, USA} \email{yau@math.harvard.edu}
\author{Xiao Zhang}
\address[Zhang]{Institute of Mathematics, Academy of Mathematics and System
Sciences, Chinese Academy of Sciences, Beijing 100080, P.R. China}
\email{xzhang@amss.ac.cn}

\begin{abstract}
We find two conditions related to the {\it news functions} of the
Bondi's radiating vacuum spacetimes. We provide a complete proof
of the positivity of the Bondi mass by using Schoen-Yau's method
under one condition and by using Witten's method under another
condition.
\end{abstract}


\date{}

\maketitle \pagenumbering{arabic}
\pagenumbering{arabic}
\mysection{Introduction}

Gravitational waves are predicted by Einstein's general
relativity. They are time dependent solutions of the Einstein
field equations which radiate or transport energy. Although they
have not been detected yet, the existence of gravitational waves
has been proved indirectly from observations of the pulsar PSR
1913+16. This rapidly rotating binary system should emit
gravitational radiation, hence lose energy and rotate faster. The
observed relative change in period
agrees remarkably with the theoretical value .\\

A fundamental conjecture is that gravitational waves can not carry
away more energy than they have initially in an isolated
gravitational system. It is usually referred as the positive mass
conjecture at null infinity. In the Bondi's radiating vacuum
spacetime, this conjecture is equivalent to the positivity of the
Bondi mass. In the pioneering work of Bondi, van der Burg, Metzner
and Sachs on the gravitational waves in vacuum spacetimes, the
Bondi mass associated to each null cone is defined and their main
result asserts that this Bondi mass is always non-increasing with
respect to the retarded time \cite{BBM, S, vdB}. Therefore, the
Bondi mass can be interpreted as the total mass of the isolated
physical system measured after the loss due to the gravitational
radiation up to that time. The proof of the positivity of the
Bondi mass was outlined by Schoen-Yau modifying their arguments in
the proof of the positivity of the ADM mass \cite{SY}. It was also
outlined by physicists applying Witten's spinor method, eg, see
\cite{IN, HP, AsHo, LV, ReuT, HT}. The main goal of this paper is
to find a complete proof of the positivity. Indeed, we find
certain conditions related to the {\it news functions} of the
system. The Bondi mass is nonnegative under these conditions.\\

It is an open problem whether vacuum Einstein field equations
always develop logarithmic singularities at null infinity. In
\cite{CMS}, the authors studied the polyhomogeneous Bondi
expansions. The $u$-evolution equations actually indicate that the
logarithmic singularities at null infinity can be removed in the
axisymmetric case (Appendix D of \cite{CMS}) if the free function
$\gamma _2 (u, x ^\alpha)$ is chosen to be zero and $\gamma _{3,1}
(u _0, x ^\alpha)$ is chosen to be zero for some $u _0$. It is
quite possible that the Bondi's radiating vacuum spacetime does
not develop any logarithmic singularity at null infinity after
suitable ``gauge fixing''. Therefore we do not consider the
polyhomogeneous Bondi expansions \cite{CJK} in the present paper.\\

The paper is organized as follows. In Section 2, we state some
well-known formulas and results of Bondi, van der Burg, Metzner
and Sachs. We employ two fundamental assumptions: {\bf Condition
A} and {\bf Condition B}. We also derive a generalized Bondi mass
loss formula under these two conditions. In Section 3, we study
some basic geometry of the asymptotically null spacelike
hypersurfaces. We compute the asymptotic behaviors of the induced
metric and the second fundamental form of an asymptotically null
spacelike hypersurface given by a certain graph. In Section 4, we
use Schoen-Yau's method to prove,  under {\bf Condition A} and
{\bf Condition B}, if there is a retarded time $u _0$ such that
$\mathcal{M}(u _0, \theta, \psi)$ defined in Section 2 is
constant, then the Bondi mass is nonnegative in the region $\{u
\leq u _0\}$, and the Bondi mass is zero at $u \in (-\infty, u
_0]$ if and only if the spacetime is flat in certain neighbourhood
of spacelike hypersurface $\{u=u _0 +\sqrt{1+r ^2}-r\}$. In
Section 5, we use Witten's method to prove, under {\bf Condition
A} and {\bf Condition B}, if there is a retarded time $u _0$ such
that $c(u _0, \theta, \psi)=d(u_0, \theta, \psi)=0$, then the
Bondi mass is nonnegative in the region $\{u \leq u _0\}$, and the
Bondi mass is zero at $u \in (-\infty, u _0]$ if and only if the
spacetime is flat in certain neighbourhood of spacelike
hypersurface $\{u=u _0 +\sqrt{1+r ^2}-r\}$. In Section 6, we
modify the definition of the Bondi energy-momentum and prove its
positivity without {\bf Condition B}.\\

\mysection{Bondi's radiating spacetimes}

We assume that $\big(\L ^{3,1}, \tilde g\big) $ is a vacuum
spacetime (possible with black holes) and $\tilde g$ takes the
following Bondi's radiating metric
 \beq
\tilde g &=&\Big(\frac{V}{r} e ^{2\beta} +r ^2 e ^{2 \gamma} U ^2
\cosh 2\delta +r ^2 e ^{-2 \gamma} W ^2 \cosh 2\delta
\nonumber\\
&&+2 r ^2 U W \sinh 2 \delta \Big)du ^2 -2e ^{2\beta}
du dr   \nonumber\\
& &-2r ^2 \Big(e ^{2 \gamma} U \cosh 2\delta +W \sinh 2 \delta
\Big) du d\theta     \nonumber\\
& &-2r ^2 \Big(e ^{-2 \gamma} W \cosh 2\delta+U \sinh 2\delta
\Big)\sin \theta du d\psi    \nonumber\\
& &+r ^2 \Big(e ^{2 \gamma} \cosh 2\delta d\theta ^2 +e ^{-2
\gamma}\cosh 2\delta \sin ^2 \theta d \psi ^2   \nonumber\\
& &+2 \sinh 2\delta \sin \theta d \theta d \psi \Big)
\label{bondi-metric}
 \eeq
in coordinates $(u, r, \theta,\psi)$ ($u$ is retarded time) where
 \beQ
-\infty < u < \infty, \;\; r >0,\;\;
 0 \leq \theta \leq \pi,\;\; 0 \leq \psi \leq 2 \pi,
 \eeQ
Denote
 \beQ
x ^0 =u,\;\;x ^1=r,\;\;x ^2=\theta, \;\;x ^3=\psi.
 \eeQ
We suppose that $\beta$, $\gamma$, $\delta$, $U$, $V$ and $W$ are
smooth functions of $u, r, \theta, \psi$. Denote $f _{,\nu} =
\frac{\partial f}{\partial x ^\nu }$ for $\nu=0,1,2,3$ throughout
the paper. The metric (\ref{bondi-metric}) was studied by Bondi,
van der Burg, Metzner and Sachs in the theory of gravitational
waves in general relativity \cite{BBM, S, vdB}. They proved that
the following asymptotic behavior holds for $r$ sufficiently large
if the spacetime satisfies the outgoing radiation condition
\cite{vdB}
 \beQ
\gamma &=&\frac{c(u, \theta, \psi)}{r}
+\frac{C(u,\theta, \psi)-\frac{1}{6} c ^3 -\frac{3}{2} c d ^2}{r ^3}
+O\big(\frac{1}{r ^4}\big),  \\
\delta &=&\frac{d(u, \theta, \psi)}{r} +\frac{H(u,\theta,
\psi)+\frac{1}{2} c ^2 d -\frac{1}{6} d ^3}{r ^3}
+O\big(\frac{1}{r ^4}\big),  \\
\beta &=&-\frac{c ^2 + d ^2}{4r ^2} +O\big(\frac{1}{r ^4}\big),     \\
U &=& -\frac{l(u, \theta, \psi)}{r ^2} +\frac{p(u, \theta,
\psi)}{r ^3}
+O\big(\frac{1}{r ^4}\big),  \\
W &=& -\frac{\bar l(u, \theta, \psi)}{r ^2}+\frac{\bar p(u,
\theta, \psi)}{r ^3}+O\big(\frac{1}{r ^4}\big),  \\
V &=& -r +2 M (u, \theta, \psi)+\frac{\bar M(u, \theta, \psi)}{r}
+O\big(\frac{1}{r^2}\big),
 \eeQ
where
 \beQ
l &=& c _{, 2} +2c \cot \theta +d _{, 3} \csc \theta,\\
\bar l &=& d _{, 2} +2d \cot \theta -c _{,3} \csc \theta,\\
p &=& 2N +3(c c _{, 2}+d d _{, 2}) +4(c ^2 +d ^2)\cot \theta\\
  & &-2(c_{,3} d -c d _{,3}) \csc \theta,\\
\bar p &=& 2P +2(c _{, 2} d -c d _{, 2}) +3(c c_{,3} +d d
_{,3})\csc \theta,\\
\bar M &=&N _{,2} +\cot \theta +P _{, 3} \csc \theta
-\frac{c ^2 +d ^2}{2}\\
& &-\big[(c _{,2}) ^2 +(d _{,2}) ^2 \big]-4(c c _{,2} +d d
_{,2})\cot \theta \\
& & -4(c ^2 +d ^2) \cot ^2
\theta -\big[(c _{,3}) ^2 +(d _{,3}) ^2 \big]\csc ^2 \theta \\
& &+4(c _{,3} d -c d _{, 3}) \csc \theta \cot \theta +2(c_{,3} d
_{,2}-c _{,2} d _{,3})\csc \theta.
 \eeQ
Here $M$ is refereed as the {\it mass aspect} and ${\it c _{,0}}$,
${\it d _{,0}}$ are refereed as the {\it news functions}. The
u-derivatives of certain functions are
 \beQ
C _{,0} &=&\frac{c ^2 c _{,0}}{2} +c d d _{,0}-\frac{c _{,0} d
^2}{2} +\frac{c M}{2} +\frac{d \lambda}{4}\\
& & -\frac{N _{,2} -N \cot\theta -P _{,3} \csc \theta}{4},\\
H _{,0} &=&-\frac{c ^2 d _{,0}}{2} +c c _{,0} d +\frac{d _{,0} d
^2}{2} +\frac{d M}{2} -\frac{c \lambda}{4}\\
& &-\frac{P _{,2} -P \cot \theta +N _{,3} \csc \theta }{4},\\
M _{,0} &=&-\big[(c _{,0} )^2 +(d _{, 0} )^2
\big]+\frac{1}{2}\big(l _{,2} +l \cot \theta +\bar l _{,3} \csc
\theta\big) _{,0}, \\
3N _{,0} &=&-M _{,2} -\frac{\lambda _{,3} \csc \theta }{2} -\big(c
_{,0} c _{,2} +d _{,0} d _{,2}\big) \\
& &-3\big(c c _{,02} +d d _{,02} \big) -4 \big(c c _{,0} +d d
_{,0} \big) \cot \theta \\
& &+\big(c _{,0} d _{,3} -c _{,3} d
_{,0}+3c _{,03} d -3 c d _{,03}\big)\csc \theta,\\
3P _{,0} &=&-M _{,3} \csc \theta +\frac{\lambda _{,2}}{2} +\big(c
_{,2} d _{,0} -c _{,0} d _{,2}\big) \\
& &+3\big(c d _{,02} - c _{,02} d \big)+4 \big(c d _{,0} -c _{,0}
d \big) \cot \theta \\
& &-\big(c _{,0} c _{,3} +d _{,0} d _{,3} +3c c _{,03} +3 d d
_{,03}\big)\csc \theta
 \eeQ
where
 \beQ
\lambda ={\bar l} _{,2} +\bar l \cot \theta -l _{,3}\csc \theta .
 \eeQ
Denote
 \beq
\mathcal{M} (u, \theta, \psi)&=& M (u, \theta,
\psi)-\frac{1}{2}\big(l _{,2} +l \cot \theta +\bar l
_{,3} \csc \theta\big) \label{mod-M}\\
 &= &M (u, \theta,
\psi) -\frac{1}{2}\big[-2c(u, \theta, \psi)+c _{,22}(u,
\theta, \psi)  \nonumber\\
& &-\csc ^2 \theta c _{,33}(u, \theta, \psi) +2 \csc \theta d
_{,23}(u,\theta, \psi) \nonumber\\
& &+3\cot \theta c _{,2} (u, \theta, \psi) +2\cot\theta \csc\theta
d_{,3}(u, \theta, \psi)\big]. \nonumber
 \eeq
Its u-derivative is
 \beq
{\mathcal{M}} _{,0} =-\big[(c _{,0} )^2 +(d _{, 0} )^2
\big].\label{Mu}
 \eeq\\

There are some physical conditions \cite{BBM, S, vdB} ensuring the
regularity of (\ref{bondi-metric}). In this paper, however, we
assume
 \begin{description}
\item[Condition A] Each of the six functions $\beta $, $\gamma$,
$\delta$, $U$, $V$, $W$ together with its derivatives up to the
second orders are equal at $\psi =0$ and $2\pi$.
\item[Condition B] For all $u$,
 \beQ
\int _0 ^{2\pi} c(u, 0, \psi) d\psi =0, \;\; \int _0 ^{2\pi} c(u,
\pi, \psi) d\psi =0. \\
 \eeQ
 \end{description}

Let $\N _{u _0}$ be a null hypersurface which is given by $u=u
_0$. The Bondi energy-momentum of $\N _{u _0}$ is defined by
\cite{BBM, CJM}:
 \beq
m _\nu (u _0) = \frac{1}{4 \pi} \int _{S ^2} M (u _0, \theta,
\psi) n ^{\nu} d S   \label{bondi-mass}
 \eeq
where $\nu =0, 1, 2, 3$, $n ^0 =1$, $n ^i $ the restriction of the
natural coordinate $x ^i$ to the unit sphere, i.e.,
 \beQ
n ^0 =1,\;\; n ^1 = \sin \theta \cos \psi,\;\; n ^2 = \sin \theta
\sin \psi,\;\; n ^3 = \cos \theta.
 \eeQ
Under {\bf Condition A} and {\bf Condition B}, we have \cite{BBM,
S, Z2}
 \beq
\frac{d}{du} m _\nu =-\frac{1}{4 \pi} \int _{S ^2} \big[(c _{,0}
)^2 +(d _{,0} )^2 \big] n ^\nu d S   \label{bondi-mass-loss}
 \eeq
for $\nu=0,1,2,3$. When $\nu=0$, this is the famous Bondi mass
loss formula.\\

The following proposition can be viewed as generalized Bondi mass
loss formula. However, it does not seem to appear in any
literature before.
 \begin{proposition} \label{bondi-e-loss}
Let $\big(\L ^{3,1}, \tilde g\big) $ be a vacuum Bondi's radiating
spacetime with metric $\tilde g$ given by (\ref{bondi-metric}).
Suppose that {\bf Condition A} and {\bf Condition B} hold. Then
 \beq
\frac{d}{du} \Big(m _0 - \sqrt{\sum _{1 \leq i \leq 3} m _i ^2 }
\Big) \leq 0. \label{bondi-energy-loss}
 \eeq
 \end{proposition}
 \pf Denote $|m| =\sqrt{m_1^2+m_2^2+m_3^2 }$. We assume $|m| \neq 0$ otherwise
it reduces to the Bondi mass-loss formula. We have
 \beQ
\frac{d}{du} \big(m _0 - |m|\big) &=& \frac{dm_0}{du}-\frac{1}{|m|}\sum _{1 \leq i \leq 3} \frac{dm_i}{du}\,m_i \\
&=& -\frac{1}{4\pi}\Big\{\int _{S ^2}\big[(c _{,0} )^2 +(d _{,0} )^2 \big] d S\\
& & -\frac{1}{|m|} \sum _{1 \leq i \leq 3} m_i\,\int _{S
^2}\big[(c _{,0} )^2 +(d _{,0} )^2 \big] n ^i d S \Big\}.
 \eeQ
Thus $\frac{d}{du} (m _0 - |m| ) \leq 0$ is equivalent to
 \beQ
\sum _{1 \leq i \leq 3} m_i\,\int _{S ^2}\big[(c _{,0} )^2 +(d
_{,0} )^2 \big] n ^i d S  \leq |m|\,\int _{S ^2}\big[(c _{,0} )^2
+(d _{,0} )^2 \big] d S.
 \eeQ
Using $(n^1)^2+(n^2)^2+(n^3)^2=1$, we obtain
 \beQ
& &\sum _{1 \leq i \leq 3}\Big\{\int _{S ^2}\big[(c _{,0} )^2 +(d
_{,0} )^2
\big] n ^i d S \Big\} ^2 \\
&\leq & \sum _{1 \leq i \leq 3} \Big\{\int _{S ^2}\big[(c _{,0}
)^2 +(d _{,0} )^2 \big] d S \Big\}\Big\{ \int _{S ^2}\big[(c
_{,0} )^2 +(d _{,0} )^2 \big] (n ^i ) ^2 d S \Big\} \\
&= &\Big\{\int _{S ^2}\big[(c _{,0} )^2 +(d _{,0} )^2 \big] d S
\Big\} ^2.
 \eeQ
Then by Cauchy-Schwarz inequality,
 \beQ
& & \sum _{1 \leq i \leq 3} m_i\,\int _{S ^2}\big[(c _{,0} )^2 +(d
_{,0} )^2
\big] n ^i d S\\
&\leq& |m| \,\sqrt{\sum _{1 \leq i \leq 3} \Big\{\int _{S
^2}\big[(c _{,0} )^2 +(d _{,0} )^2
\big] n ^i d S\Big\} ^2}\\
&\leq& |m|\,\int _{S ^2}\big[(c _{,0} )^2 +(d _{,0} )^2 \big] d S.
 \eeQ
Therefore (\ref{bondi-energy-loss}) holds. \qed\\
\mysection{Asymptotically null spacelike hypersurfaces}

The hypersurface $u=\sqrt{1+r ^2} -r $ in the Minkowski spacetime
is a hyperbola equipped with the standard hyperbolic 3-metric
$\breve{g}$. Let $\{\breve{e} _i \}$ be the frame
 \beQ
\breve{e} _1 = \sqrt{1+r ^2} \frac{\pa}{\partial r},
\quad\breve{e} _2 = \frac{1}{r} \frac{\pa}{\partial \theta}, \quad
\breve{e} _3 = \frac{1}{r \sin \theta} \frac{\pa}{\partial \psi}.
 \eeQ
Let $\{ \breve{e} ^i \}$ be the coframe. Denote $\breve{\nabla }
_{i} =\breve{\nabla } _{\breve{e} _i}$, etc., where
$\breve{\nabla} $ is the Levi-Civita connection of $\breve{g}$.
The connection 1-forms $\{\breve{\omega} _{ij}\}$ are given by $d
\breve{e} ^{i}= - \breve{\omega} _{ij} \wedge \breve{e} ^{j}$ or
$\breve{\nabla } \breve{e} _{i}= - \breve{\omega} _{ij} \otimes
\breve{e} _{j}$. They are
 \beQ
\breve{\omega}  _{12} = -\frac{\sqrt{1+r ^2}}{r}\;\breve{e}
^{2},\quad \breve{\omega}  _{13} = -\frac{\sqrt{1+r
^2}}{r}\;\breve{e}^{3},\quad \breve{\omega}  _{23} =
-\frac{\cot\theta}{r}\;\breve{e}^{3}.
 \eeQ

Let $\X$ be a spacelike hypersurface in vacuum Bondi's radiating
spacetime $\big(\L ^{3,1}, \tilde g \big)$ with metric
(\ref{bondi-metric}), which is given by the inclusion:
 \beQ
\i: \X ^3  & \longrightarrow & {\L} ^{3,1}\\
         \big(y ^1,y^2, y^3 \big)&\longmapsto &\big(x ^0, x ^1, x^2, x^3 \big)
 \eeQ
for $r$ sufficiently large, where
 \beQ
 x ^0 = u(y ^1, y ^2, y  ^3),\;\;
 x ^1 = y ^1 =r,\;\;
 x ^2 = y ^2 =\theta,\;\;
 x ^3 = y ^3 =\psi.
 \eeQ
Let $g =\i ^{*} \tilde g$ be the induced metric of $\X$ and $h$ be
the second fundamental form of $\X$. Let $\tilde \nabla$ be the
Levi-Civita connection of $\L ^{3,1}$. For any tangent vectors $Y
_i, Y _j \in T\X$, $\i _{*} Y _i$, $\i _{*} Y _j$ are the tangent
vectors along $\X$,
 \beQ
g\big(Y_i, Y _j\big)=\tilde g \big(\i _{*} Y _i, \i _{*} Y _j
\big).
 \eeQ
Let $e _n$ be the downward unit normal of $\X$. The second
fundamental form is defined as
 \beQ
h\big(Y_i, Y _j\big)=\tilde g \big(\tilde{\nabla } _{\i _{*} Y _i}
\i _{*} Y _j, e _n \big).
 \eeQ
Now it is a straightforward computation that
 \beQ
\i _{*} \frac{\pa}{\pa y ^i} =\frac{\pa}{\pa x ^\alpha}\frac{\pa x
^\alpha}{\pa y ^i} = \frac{\pa}{\pa x ^0}\frac{\pa x
^0}{\pa y ^i}+\frac{\pa}{\pa x ^i}.
 \eeQ
Denote $e _i =\i _{*} \breve{e} _i$. Then
 \beq
e _1 =\sqrt{1+r ^2} u _{,1}\, \pa _0 +\breve{e} _1,\quad e _2
=\frac{u _{,2}}{r} \,\pa _0  +\breve{e} _2, \quad e _3 =\frac{u
_{,3}}{r\sin \theta}\, \pa _0 +\breve{e} _3. \label{e}
 \eeq

\begin{definition}
A spacelike hypersurface $(\X, g, h)$ in an asymptotically flat
spacetime is asymptotically null of order $\tau >0$ if, for $r$
sufficiently large, $g(\breve{e} _i,\breve{e} _j)=\delta
_{ij}+a_{ij}$, $h(\breve{e} _i, \breve{e} _j)=\delta
_{ij}+b_{ij}$, where $a _{ij}$, $b _{ij}$ satisfy
 \beq
\big\{a_{ij},\breve{\nabla} _k a_{ij}, \breve{\nabla} _l
\breve{\nabla} _k a_{ij},
 b_{ij}, \breve{\nabla} _k
 b_{ij}\big\}=O\big(\frac{1}{r^\tau}\big).
\label{a-b}
 \eeq
\end{definition}

Let $(\X, g, h)$ be an asymptotically null spacelike hypersurface
with the induced metric $g$ and the second fundamental form $h$ in
vacuum Bondi's radiating spacetime $\big(\L ^{3,1}, \tilde g
\big)$, which is given by
 \beq
u =\sqrt{1+r ^2} -r + \frac{\big(c ^2 +d ^2\big) _{u=0}}{12r ^3}
+\frac{a _{3} (\theta, \psi)}{r ^4}+a _{4} \label{u-surface-0}
 \eeq
where $a _{4} (r,\theta, \psi)$ is a smooth function which
satisfies: In the Euclidean coordinate systems $\{\breve{z} ^i\}$,
$|\breve{z}| =r$,
 \beQ
a_4 =o\big(\frac{1}{r ^4}\big), \quad\pa _k a_4 =o\big(\frac{1}{r
^5}\big), \quad\pa _k\pa _l a_4 =o\big(\frac{1}{r ^6}\big)
 \eeQ
as $r \rightarrow \infty$. We will compute asymptotic behaviors of
the induced metric and the second fundamental form of $\X$. The
induced metric can be obtained by substituting $du$ into
(\ref{bondi-metric}). Let $X _n$ be the downward normal vector
 \beQ
X _n = -\frac{\pa }{\pa x ^0}- \varrho ^i \frac{\pa }{\pa x ^i}.
 \eeQ
Let $e_i$ be given by (\ref{e}). Since $X _n$ is orthogonal to $e
_i$, we obtain
 \beQ
\tilde g (e _i, X_n) =0.
 \eeQ
This implies that $\varrho ^i$ satisfies the following equations
 \beQ
\big(u _{,i}\; \tilde g_{00} +\tilde g _{0i} \big)
    +\varrho ^j \big(u _{,i} \;\tilde g_{0j} +\tilde g _{ij}\big)=0
 \eeQ
for $i=1,2,3$. Therefore $\varrho ^i$ can be found by solving
linear algebraic equations and the unit normal vector is
 \beQ
e _n = \frac{X_n}{\sqrt{-\tilde g(X _n, X_n)}}.
 \eeQ
The second fundamental form is then given by
  \beQ
h (\breve{e} _i, \breve{e} _j) = \tilde g \big(\tilde{\nabla } _{e
_i} e _j, e _n \big)
 \eeQ
for $1 \leq i, j \leq 3$. Now we define $a \approx b$ if and only
if $a =b +o\big(\frac{1}{r^3}\big)$. For $r$ sufficiently large,
we expand $c$, $d$ and $M$ at $u=0$ by Taylor series
 \beq
c(u, \theta, \psi)&\approx &c(0, \theta, \psi)+c _{,0}(0, \theta,
\psi) u +\frac{c
_{,00} (0, \theta, \psi)}{2}  u ^2,  \label{c}\\
d(u, \theta, \psi)&\approx &d (0, \theta, \psi)+d _{,0} (0,
\theta, \psi) u +\frac{d
_{,00}(0, \theta, \psi)}{2} u ^2,  \label{d}\\
M(u, \theta, \psi)&\approx&M (0, \theta, \psi)+M _{,0} (0, \theta,
\psi) u +\frac{M _{,00}(0, \theta, \psi)}{2}  u ^2, \label{M}
 \eeq
with the help of Mathematica 5.0, we obtain the asymptotic
behaviors of the metric $g$
 \beQ
g(\breve{e} _1, \breve{e} _1) & \approx &1+\frac{16a_3+M -c
c_{,0}- d d
_{,0}}{2r^3},\\
 g(\breve{e} _1, \breve{e} _2)
&\approx &-\frac{l}{2r^2}
 +\frac{12 N -3 l _{,0}
 +4 (c c _{,2}+ d d _{,2})}{12r ^3},\\
g(\breve{e} _1, \breve{e} _3)&\approx&-\frac{\bar l}{2r ^2}
+\frac{12 P-3\bar{l} _{,0}+4\csc \theta (c c _{,3}+  d d _{,3})}{12r ^3},\\
g(\breve{e} _2, \breve{e} _2)&\approx&1+\frac{2c}{r}+\frac{2(c^2+ d^2) +c _{,0}}{r ^2}\\
  & & +\frac{c ^3 +c d ^2 +2C+2(c c _{,0} + d d_{,0})
         +\frac{c _{,00}}{4}}{r ^3},\\
g(\breve{e} _2, \breve{e} _3) &\approx
&\frac{2d}{r}+\frac{d_{,0}}{r^2}
   +\frac{c ^2 d +d ^3 +2H+\frac{d _{,00}}{4}}{r ^3}\\
g(\breve{e} _3, \breve{e} _3) &\approx &1-\frac{2c}{r}+\frac{2(c^2
+d^2) -c _{,0}}{r ^2}\\
     & &+\frac{-c ^3 -c d ^2-2C +2(c c _{,0} +d d _{,0}) -\frac{c
_{,00}}{4}}{r ^3}.\\
 h(\breve{e} _1,
\breve{e} _1)
   &\approx&1 + \frac{c ^2+d^2}{r ^2}
   +\frac{16a_3-M}{r ^3},\\
h(\breve{e} _1, \breve{e} _2)
   &\approx & \frac{l}{2r ^2}
      +\frac{1}{2r^3}\big[ \frac{l _{,0}}{2}-2 (c ^2+ d ^2)  \cot\theta -4 N\\
   & & (-c d _{,3} +c _{,3}d )\csc \theta -\frac{13}{3}(c c _{,2} +d d _{,2})\big],\\
 h(\breve{e} _1, \breve{e} _3)
  &\approx & \frac{\bar l}{2r^2 }
      +\frac{1}{2r ^3 }\big[\frac{\bar l _{,0}}{2} +c d _{,2}-c _{,2} d
      -4P\\
   & & -\frac{13}{3} (c c _{,3} +d d _{,3}) \csc \theta\big],\\
h(\breve{e} _2, \breve{e} _2)&\approx &1 +\frac{c}{r}
     +\frac{c _{,0}}{r ^2}
   +\frac{1}{4r ^3}\big[ 3M-16 a_3 -4C -2 l _{,2} \\
   & &- 2c(c ^2 + d ^2) +5(c c _{,0} +d d _{,0})+\frac{3}{2}c
   _{,00}\big],\\
h(\breve{e} _2, \breve{e} _3)
       &\approx&\frac{d}{r}+\frac{d _{,0}}{r^2}
        +\frac{1}{4r ^3 }\big[-2d(c ^2 +d ^2) +2d \cot ^2 \theta
        \\
& &+2d \csc ^2 \theta -4c _{,3} \cot \theta \csc \theta -d
_{,33}\csc ^2 \theta \\
& & -d _{,2}\cot \theta - d _{,22}
          - 4H +\frac{3}{2} d_{,00} \big],\\
h(\breve{e} _3, \breve{e} _3)&\approx &1 -\frac{c}{r}
     -\frac{c _{,0}}{r ^2}
   +\frac{1}{4r ^3}\big[ 3M-16a_3 +4C \\
   & &+2c(c^2 + d ^2) +5(cc _{,0}
   +dd _{,0})-\frac{3}{2}c _{,00}\\
   & &-2 l \cot \theta -2 \bar l _{,3} \csc \theta \big].
 \eeQ
Here all functions in the right hand sides take value at $u=0$ and
all derivatives with respect to $x ^2$ and $x ^3$ are taken after
substituting $u=0$. Therefore $(\X, g, h)$ is asymptotically null of order $1$.\\
\mysection{Positivity - Schoen-Yau's method}

In this section, we will complete the argument in \cite{SY}.
Denote by $(\X, g, h)$ the asymptotically null spacelike
hypersurface which is given by (\ref{u-surface-0}) for $r$
sufficiently large. In \cite{SY}, Schoen-Yau solved the following
Jang's equation on $\X$:
 \beq
\Big(g ^{ij} -\frac{f^i  f^j}{1+|\nabla f | ^2}\Big)\Big(\frac{f
_{,ij}}{\sqrt{1+|\nabla f | ^2}} -h _{ij} \Big) =0 \label{jang-eq}
 \eeq
under the suitable boundary condition
 \beq
f \rightarrow f _0 \label{j-f0}
 \eeq
as $r \rightarrow \infty$ such that the metric
 \beq
\bar g = g +\nabla f \otimes \nabla f   \label{bar-g}
 \eeq
is asymptotically flat. Denote by $J(f)$ the left hand side of
Jang's equation (\ref{jang-eq}). Note that in the standard
hyperbolic 3-space, (\ref{jang-eq}) has a solution $f=\sqrt{1+r
^2}$. Therefore it is reasonable to set
 \beQ
f _0 =\sqrt{1+r ^2} +o(r).
 \eeQ
Let $f$ be a function on $\X$ which has asymptotic expansion
 \beq
f=\sqrt{1+r ^2}+p(\theta, \psi)\ln r +q (r, \theta, \psi),
\label{f}
 \eeq
for $r$ sufficiently large, where $p(\theta, \psi)$ is a smooth
function on $S ^2$ and $q$ is a smooth function on $\R ^3$ which
satisfies the following asymptotic conditions: In the Euclidean
coordinate systems $\{\breve{z} ^i\}$, $|\breve{z}| =r$,
 \beQ
q=o\big(1\big),\quad \pa _k q =o\big(\frac{1}{r}\big), \quad\pa _k
\pa _l q =o\big(\frac{1}{r ^2}\big),\quad \pa _k \pa _l\pa _j q
=o\big(\frac{1}{r ^3}\big)
 \label{asy-f}
 \eeQ
as $r \rightarrow \infty$.\\

Let the standard metric of $S ^2$ be $d\theta ^2 +\sin ^2 \theta
d\psi ^2$. The Laplacian operator of this metric is
 \beQ
\triangle _{S ^2} =\frac{\pa ^2}{{\pa \theta} ^2} +\cot \theta
\frac{\pa }{\pa \theta} +\csc ^2 \theta \frac{\pa ^2}{{\pa \psi}
^2}.
 \eeQ
The {\it spherical harmonics $w _j$} are the eigenfunctions of
$\triangle _{S ^2}$, i.e., $\triangle _{S ^2} w _j =j(j-1) w_j$
for $j=1,2,\cdots$.\\
\begin{proposition}
If Jang's equation (\ref{jang-eq}) has a solution $f$ which has
the asymptotic expansion (\ref{f}) for $r$ sufficiently large,
then $p(\theta, \psi)$ and ${\mathcal M}(0, \theta, \psi)$ must be
constant.
\end{proposition}
\pf A lengthy computation with the help of Mathematica 5.0 shows
that
 \beQ
J(f)\approx \frac{\ln r}{r ^3} \triangle _{S ^2}
p+\frac{p-2\mathcal{M}(0, \theta, \psi)}{r ^3}
 \eeQ
for $r$ sufficiently large. That $J(f)=0$ implies
 \beQ
\triangle _{S ^2} p=0, \quad p-2 \mathcal{M}=0.
 \eeQ
As there is no nonconstant harmonic function on $S ^2$, the
proposition follows. \qed\\

The existence of (\ref{jang-eq}) under the boundary condition
(\ref{j-f0}) with
 \beq
f _0 (r)=\sqrt{1+r ^2} +p \ln r \label{f0}
 \eeq
for certain constant $p$ can be established as follows: We extend
$f _0$ to the whole $\X$ and denote as $f _0$ also. Denote $B _R$
as the ball of radius $R$ in $\R ^3$. If $\X$ has no apparent
horizon, the existence theorem for the Dirichlet problem \cite{Y}
indicates that there exists a (smooth) solution $\bar f _R$ of
(\ref{jang-eq}) on $B _R$ such that
 \beQ
\bar f _R   \big | _{\partial B _R} =0
 \eeQ
for sufficiently large $R$. By the translation invariance of
(\ref{jang-eq}) in the vertical direction, we find that
 \beQ
f _R = \bar f _R +f _0 (R)
 \eeQ
is a solution of (\ref{jang-eq}) which is $f _0 (R)$ on $\partial
B _R$. Now the estimates in \cite{SY0} show that
 \beQ
f _R \longrightarrow f
 \eeQ
on any compact subset of $\X$, where $f$ is a (smooth) solution of
(\ref{jang-eq}). Write $f=f _0 +f _1$ where $\lim _{r \rightarrow
\infty} f _1 =0$. Substitute it into Jang's equation
(\ref{jang-eq}) and obtain an equation for $f _1$. Then use the
similar argument as the proof of Proposition 3 in \cite{SY0}, we
can show that for any $\varepsilon \in (0,1)$, there is a constant
$C(\varepsilon)$ depending only on $\varepsilon $ and the geometry
of $\X$ such that
 \beQ
\big|f_1 (\breve{z})\big| + |\breve{z}| \big|\pa f_1
(\breve{z})\big| +|\breve{z}| ^2\big|\pa \pa f_1
(\breve{z})\big|+|\breve{z}| ^3 \big|\pa \pa \pa f_1
(\breve{z})\big| \leq C(\varepsilon) |\breve{z}| ^\varepsilon .
 \eeQ
Therefore $f$ has asymptotic behaviors (\ref{f}), (\ref{asy-f})
for $r$ sufficiently large.

By adding one point compactification, the existence of
(\ref{jang-eq}) can be extended to $\X$ with apparent horizons.
See \cite{SY0} for detail. \\

The following lemma was proved in \cite{Z2}.
\begin{lemma}\label{L}
Let $\big(\L ^{3,1}, \tilde g\big) $ be a vacuum Bondi's radiating
spacetime with metric $\tilde g$ given by (\ref{bondi-metric}).
Suppose that {\bf Condition A} and {\bf Condition B} hold. Then
 \beQ
\int _{S ^2}\big(l _{,2} +l \cot \theta +\bar l _{,3} \csc \theta
\big) n ^\nu  dS =0
 \eeQ
for $\nu =0,1,2,3$.\\
\end{lemma}

 \begin{theorem}\label{positivity-bondi-mass-2}
Let $\big(\L ^{3,1}, \tilde g\big) $ be a vacuum Bondi's radiating
spacetime with metric $\tilde g$ given by (\ref{bondi-metric}).
Suppose that {\bf Condition A} and {\bf Condition B} hold. If
there exists a constant $u _0$ such that $\mathcal{M} (u _0,
\theta, \psi)$ is constant, then
 \beQ
m _{0} (u) \geq \sqrt{\sum _{1 \leq i \leq 3} m ^2 _{i} (u)}
 \eeQ
for all $u \leq u_0$. If the equality holds for some $u \in
(-\infty, u _0]$, $\L ^{3,1}$ is flat in the region foliated by
all spacelike hypersurfaces which are given by
 \beQ
u=u _0 + \sqrt{1+r ^2} -r +o(\frac{1}{r^4})
 \eeQ
for $r$ sufficiently large. In particular, if the equality holds
for all $u \leq u _0$, $\L ^{3,1}$ is flat in the region $\{u \leq
u _0\}$.
 \end{theorem}
\pf Suppose $\mathcal{M} (u _0, \theta, \psi)=\frac{p}{2}$. By the
translation invariance of Jang's equation, we can assume that $u_0
=0$. The assumption of the theorem ensures that there exists a
smooth solution $f$ of Jang's equation (\ref{jang-eq}) under the
boundary condition (\ref{j-f0}) with $f _0$ given by (\ref{f0}).
It is obvious that the metric $\bar g$ given by (\ref{bar-g}) is
asymptotic flat. Now we show its ADM total energy is $p$. Denote
by $g _0$ the flat metric of $\R ^3$ in polar coordinates. Let
$\{e ^0 _i \}$ be the frame of $g _0$
 \beQ
 e ^0 _1 = \frac{\pa}{\partial r},\quad
 e ^0 _2 = \frac{1}{r} \frac{\pa}{\partial \theta},\quad
 e ^0 _3 = \frac{1}{r \sin \theta} \frac{\pa}{\partial \psi}.
 \eeQ
Let $\{ e_0 ^i \}$ be the coframe of $g _0$. Denote $\alpha
_{ij}=\bar g \big(e ^0 _i, e ^0 _j \big)-\delta _{ij}$. Now we use
the ADM energy expression in polar coordinates
 \beQ
E(\bar g)&=&\frac{1}{16\pi}\lim _{r \rightarrow \infty} \int _{S
_r} \big[{(\nabla ^0 )} ^j \alpha _{1j} -(\nabla ^0 ) _1 tr _{g
_0} ( \alpha ) \big] {e _0} ^2 \wedge {e _0} ^3
 \eeQ
where $\nabla ^0$ is the Levi-Civita connection of $g _0$. Since
 \beQ
{(\nabla ^0 )} ^j \alpha _{1j} -(\nabla ^0 ) _1 tr _{g _0} (
\alpha )=\frac{\ln r}{r ^2} \triangle _{S ^2} p+\frac{4p}{r ^2}
+o\big(\frac{1}{r^2}\big),
 \eeQ
we obtain
 \beQ
E(\bar g)=p.
 \eeQ
Since it satisfies vacuum Einstein field equations, the Bondi's
radiating metric satisfies the dominant energy condition
automatically. Therefore the scalar curvature $\bar R$ of $\bar g$
satisfies
 \beQ
\bar R \geq 2 \big| Y \big| ^2 _{\bar g} -2 div _{\bar g} Y
 \eeQ
for certain vector field in $\bar \X$. Therefore a standard
positive mass argument \cite{SY0, LY} shows that
 \beQ
E(\bar g) =p \geq 0.
 \eeQ
And $p=0$ if and only if the metric $\bar g$ is flat which implies
that $(\X, g, h)$ can be embedded into the Minkowski spacetime as
a spacelike hypersurface with the induced metric $g$ from the
Minkowski metric and the second fundamental form $h$.

Integrating $\mathcal{M}(0,\theta, \psi)=\frac{p}{2}$ over unit $S
^2$ and using Lemma \ref{L}, we obtain the Bondi energy-momentum
of slice $u=0$
 \beQ
m _0 (0) =\frac{p}{2},\;\;m _1 (0) =m _2 (0)=m _{3}(0)=0.
 \eeQ
Thus the theorem follows from Proposition \ref{bondi-e-loss}.
\qed\\

\mysection{Positivity - Witten's method}

In this section, we will use Witten's method \cite{W} and the
positive mass theorem near null infinity proved by the third author
\cite{Z1, Z3} to study the positivity of the Bondi mass. Let
$\big(\X, g, h\big)$ be an asymptotically null spacelike
hypersurface. Denote
 \beQ
E_{\nu } (\X)&=& \frac{1}{16\pi}\lim _{r \rw \infty} \int _{S_r}
\E n ^\nu r \;\breve{e} ^2 \wedge \breve{e} ^3,
\label{hyperbolic-mass}\\
P_{\nu } (\X)&=& \frac{1}{8\pi}\lim _{r \rw \infty} \int _{S_r} \P
n ^\nu r \;\breve{e} ^2 \wedge \breve{e} ^3
\label{hyperbolic-momentum}
 \eeQ
where
 \beQ
 \E &=&\breve{\nabla} ^ {j} a_{1j} - \breve{\nabla} _1 tr _{\breve{g}} (a)
-\big[a_{11}-\delta _{11} tr _{\breve{g}} (a) \big], \\
 \P &=& b_{11}-\delta _{11} tr _{\breve{g}}(b).
 \eeQ
Theorem 4.1 in \cite{Z1} indicates if $(\X, g, h)$ is
asymptotically null spacelike hypersurface of order $\tau
>\frac{3}{2}$ in vacuum Bondi's radiating spacetime (\ref{bondi-metric}),
then
 \beq
E _0 (\X) - P _0 (\X) \geq \sqrt{\sum _{1 \leq i \leq 3} \big[E _i
(\X) - P _i (\X)\big] ^2} \label{positive-mass}
 \eeq
and the equality implies the spacetime is flat over $\X$. (Theorem
4.1 was proved for $\tau =3$. However, the argument goes through
if $c | _{u=0} =d | _{u=0} =0$ for the above $(\X, g, h)$ in the
Bondi's radiating spacetimes. See also Theorem 3.1 and Remark 3.1
in \cite{Z3}. The sharp order $\tau >\frac{3}{2}$ together with
certain integrable conditions was also given in \cite{CN, CH} to
ensure the argument to work.) In general, the hyperbolic mass of
an asymptotically null spacelike hypersurface is different from
the Bondi mass of the null cone. For instance, if $c | _{u=0}$ or
$d | _{u=0}$ is nonzero, $E _0 (\X) - P _0 (\X)$ may not be
finite.\\

\begin{lemma}\label{E-P}
Let $\big(\L ^{3,1}, \tilde g\big) $ be a vacuum Bondi's radiating
spacetime with metric $\tilde g$ given by (\ref{bondi-metric}).
Let $\big(\X, g, h\big)$ be a spacelike hypersurface $u$ which is
given by (\ref{u-surface-0}) for $r$ sufficiently large. Denote
$L(\phi, \psi)=l(0,\phi, \psi)$, $\bar{L}(\phi,
\psi)=\bar{l}(0,\phi, \psi)$. Then
 \beQ
 \E &\approx&\frac{12}{r^2}\big(c ^2 + d^2\big) _{u=0}\\
    & &+\frac{1}{r ^3}\big(M+16 a_3+15c c _0 +15 d d _0 \big)_{u=0}\\
    & &-\frac{1}{2 r ^3}\big(L _{,2} +L \cot \theta +\bar L _{,3} \csc\theta
    \big),\\
 \P &\approx&-\frac{1}{2 r ^3}\big(3M -16 a_3+5c c _0 +5 d d _0
 \big)_{u=0}\\
    & &+\frac{1}{2 r ^3}\big(L _{,2} +L \cot \theta +\bar L _{,3} \csc\theta
    \big).
 \eeQ
\end{lemma}
\pf Note that
 \beQ
a_{22}+a_{33}&\approx&\frac{4}{r ^2} \big(c ^2 +d ^2\big)_{u=0}
             +\frac{4}{r ^3} \big(c c_{,0} +d d _{,0}\big)_{u=0},\\
b_{22}+b_{33}&\approx&\frac{1}{2r ^3}\big(3M -16a _3 +5c c_{,0}
+5d d
_{,0}\big)_{u=0}\\
& &-\frac{1}{2r ^3}\big(L _{,2} +L \cot \theta +\bar L _{,3} \csc
\theta \big).
 \eeQ
Using the formula
 \beQ
 \breve{\nabla} _k a _{ij} = \breve{e} _k (a _{ij})
 -a _{jl} \breve{\omega} _{li} (\breve{e} _k)
 -a _{il} \breve{\omega} _{lj} (\breve{e} _k),
 \eeQ
we obtain
 \beQ
\E &=&\breve{e} _j (a _{1j})
 -a _{jl} \breve{\omega} _{l1} (\breve{e} _j)
 -a _{1l} \breve{\omega} _{lj} (\breve{e} _j)
 -\breve{\nabla} _1 tr _{\breve{g}} (a)\\
& & +a _{22} +a _{33}\\
&\approx&-\frac{1}{2r ^3}\big(L _{,2} +L \cot \theta +\bar L _{,3}
\csc
\theta \big)\\
& &+\frac{1}{r ^3}\big(M+16 a_3 -c c _{,0} -d d _{,0}
\big)_{u=0}\\
& &+\frac{8\sqrt{1+r ^2}}{r ^3}\big(c ^2 +d ^2\big)_{u=0}\\
& &+\frac{12\sqrt{1+r ^2}}{r ^4}\big(c c_{,0} +d d_{,0}\big)_{u=0}\\
& &+\frac{4}{r ^3}\big(c ^2 +d ^2\big)_{u=0} +\frac{4}{r ^4}\big(c
c_{,0} +d d_{,0}\big)_{u=0}\\
&\approx&-\frac{1}{2r ^3}\big(L _{,2} +L \cot \theta +\bar L _{,3}
\csc
\theta \big)\\
& &+\frac{1}{r ^3}\big(M+16 a_3 +15 c c _{,0} +15 d d _{,0} \big)_{u=0}\\
& &+\frac{12}{r ^2}\big(c ^2 +d ^2\big)_{u=0}+O\big(\frac{1}{r^4}\big),\\
\P &=&-b _{22} -b _{33}\\
&\approx&-\frac{1}{2r ^3}\big(3M -16a _3 +5c c_{,0} +5d d
_{,0}\big)_{u=0}\\
& &+\frac{1}{2r ^3}\big(L _{,2} +L \cot \theta +\bar L _{,3} \csc
\theta \big).
 \eeQ
\qed

 \begin{theorem}\label{positivity-bondi-mass}
Let $\big(\L ^{3,1}, \tilde g\big) $ be a vacuum Bondi's radiating
spacetime with metric $\tilde g$ given by (\ref{bondi-metric}).
Suppose that {\bf Condition A} and {\bf Condition B} hold and $c |
_{u=u _0} = d | _{u=u _0} = 0$ for some $u _0$. Then
 \beQ
m _{0} (u) \geq \sqrt{\sum _{1 \leq i \leq 3} m ^2 _{i} (u)}
 \eeQ
for all $u \leq u _0$. If the equality holds for some $u \in
(-\infty, u _0]$, $\L ^{3,1}$ is flat in the region foliated by
all spacelike hypersurfaces which are given by
 \beQ
u=u _0 + \sqrt{1+r ^2} -r +o(\frac{1}{r^4})
 \eeQ
for $r$ sufficiently large. In particular, if the equality holds
for all $u \leq u _0$, $\L ^{3,1}$ is flat in the region $\{u \leq
u _0\}$.
 \end{theorem}
 \pf By translation, we can assume that $u _0 =0$. Choose an
asymptotically null spacelike hypersurface $\X $ which is given by
(\ref{u-surface-0}) with $a _3=0$ for $r$ sufficiently large. By
Lemma \ref{E-P}, we obtain
 \beQ
\E - \P \approx -\frac{1}{r ^3}\big(L _{,2} +L \cot \theta +\bar L
_{,3} \csc \theta \big)+\frac{5M(0,\theta,\psi)}{2r ^3}.
 \eeQ
Then Lemma \ref{L} implies that
 \beQ
E _{\nu} (\X)-P _{\nu} (\X) =\frac{5}{8} m _\nu (0).
 \eeQ
Therefore the first part of the theorem follows from
(\ref{positive-mass}) and Proposition \ref{bondi-e-loss}. For the
second part, if the equality holds for some $u \in (-\infty, u
_0]$, then the equality holds for $u =u _0$ by Proposition
\ref{bondi-e-loss}. Thus $\L ^{3,1}$ is flat over $\X$ and it follows. \qed\\
\mysection{Modified Bondi energy-momentum}

We can modify the definition of the Bondi energy-momentum to
remove {\bf Condition B}. Define the modified Bondi
energy-momentum as
 \beq
{\bf m} _\nu (u _0) = \frac{1}{4 \pi} \int _{S ^2} \mathcal{M} (u
_0, \theta, \psi) n ^{\nu} d S \label{m-bondi-mass}
 \eeq
for $\nu =0, 1, 2, 3$. Then we can prove that
 \beq
\frac{d}{du} {\bf m} _\nu =-\frac{1}{4 \pi} \int _{S ^2} \Big((c
_{,0} )^2 +(d _{,0} )^2 \Big) n ^\nu d S \label{m-bondi-mass-loss}
 \eeq
for $\nu=0,1,2,3$, and
 \beq
\frac{d}{du} \Big({\bf m} _0 - \sqrt{\sum _{1 \leq i \leq 3} {\bf
m} _i ^2 } \Big) \leq 0\label{m-bondi-energy-loss}
 \eeq
under {\bf Condition A} only.\\

Now choosing the spacelike asymptotically null hypersurface $\X $
given by (\ref{u-surface-0}) with $a _3 =-\frac{M(0, \theta,
\psi)}{16}$. If $c | _{u=0} = 0$, $d | _{u=0} = 0$, then
 \beQ
\E - \P \approx \frac{2{\mathcal {M}}(0, \theta, \psi)}{r ^3}.
 \eeQ
Therefore the following theorems are a direct consequence.\\
 \begin{theorem}
Let $\big(\L ^{3,1}, \tilde g\big) $ be a vacuum Bondi's radiating
spacetime with metric $\tilde g$ given by (\ref{bondi-metric}).
Suppose that {\bf Condition A} holds. If (i) either ${\mathcal
M}(u _0, \theta, \psi)$ is constant, (ii) or $c | _{u=u _0} = d |
_{u=u _0} = 0$ for some $u_0$, then
 \beQ
{\bf m} _{0} (u) \geq \sqrt{\sum _{1 \leq i \leq 3} {\bf m} ^2
_{i} (u)}
 \eeQ
for all $u \leq u _0$. If the equality holds for some $u \in
(-\infty, u _0]$, $\L ^{3,1}$ is flat in the region foliated by
all spacelike hypersurfaces which are given by
 \beQ
u=u _0 + \sqrt{1+r ^2} -r -\frac{M (u _0, \theta,\psi)}{16 r
^4}+o(\frac{1}{r^4})
 \eeQ
for $r$ sufficiently large. In particular, if the equality holds
for all $u \leq u _0$, $\L ^{3,1}$ is flat in the region $\{u \leq
u _0 -\frac{M (u _0, \theta,\psi)}{16 r ^4}\}$.
 \end{theorem}

\hspace{8mm}

{\footnotesize {\it Acknowledgements.} {Xiao Zhang is partially
supported by National Natural Science Foundation of China under
grants 10231050, 10421001 and the Innovation Project of Chinese
Academy of Sciences. This work was partially done when Wen-ling
Huang visited the Morningside Center of Mathematics, Chinese Academy
of Sciences, and she would like to thank the center for its
hospitality. Part of the main results was announced in \cite{Z3}.}}

\hspace{8mm}

\end{document}